\documentclass[12pt]{article}

\usepackage{amsmath}
\usepackage{amssymb}
\usepackage{amsfonts}
\usepackage{amsthm}
\usepackage{latexsym}
\usepackage{epsfig}
\usepackage{array}
\usepackage{boxedminipage}
\makeatletter
\renewcommand{\section}{\@startsection
{section} {1} {0mm} {-\baselineskip} {0.5\baselineskip}
{\large\bf}}
\renewcommand{\subsection}{\@startsection
{subsection} {2} {0mm} {-\baselineskip} {0.5\baselineskip}
{\normalsize\bf}} \makeatother

\newcommand{\R}{\mathbb{R}}
\newtheorem{defin}{Definition}
\newtheorem{teorema}{Theorem}
\newtheorem{lema}{Lemma}

\newcommand{\ke}{\mathrm{ker} \hspace{0.5mm}}
\newcommand{\im}{\mathrm{im} \hspace{0.5mm}}

\begin{document}

\mbox{}

\begin{center}
{\Large\bf  Dynamical properties of electrical circuits \vspace{2mm}\\
with fully nonlinear memristors\footnote{Supported by 
Research Project MTM2007-62064 (MEC, Spain).}} 

\ \vspace{1mm} \\

\ \\

{\large\sc Ricardo Riaza}\\ 
\ \vspace{-3mm} \\
Departamento de Matem\'{a}tica Aplicada a las Tecnolog\'{\i}as de la
Informaci\'{o}n \\
Escuela T\'{e}cnica Superior de Ingenieros de Telecomunicaci\'{o}n \\
Universidad Polit\'{e}cnica de Madrid  
- 28040 Madrid, Spain \\
{\sl ricardo.riaza@upm.es} \\

\vspace{0.5cm}


\


\end{center}

\begin{abstract}
The recent design of a nanoscale device with a memristive characteristic
has had a great impact in nonlinear circuit theory. Such a
device, whose existence was predicted by Leon Chua in 1971,
is governed by a charge-dependent voltage-current relation 
of the form $v=M(q)i$.
In this paper we show that 
allowing for a fully nonlinear characteristic $v=\eta(q, i)$ in memristive devices
provides a general framework for modeling and analyzing a very broad family
of electrical and electronic circuits; Chua's memristors 
are particular instances in which $\eta(q,i)$ is linear in $i$. 
We examine several
dynamical features of circuits with fully nonlinear memristors, accommodating
not only charge-controlled but also flux-controlled ones, with
a characteristic of the form $i=\zeta(\varphi, v)$. Our results apply in
particular to Chua's memristive circuits; certain properties
of these can be seen as a consequence of the special form
of the elastance and
reluctance matrices displayed by Chua's memristors.
\end{abstract}

\vspace{5mm}


\vspace{5mm}









\vspace{0.5cm}

\newpage

\section{Introduction}
\label{sec-intr}

In 1971 Leon Chua predicted the existence of a fourth basic circuit element
which would be governed by a flux-charge relation, having either a charge-controlled
form  $\varphi = \phi(q)$ or a flux-controlled one 
$q = \sigma(\varphi)$ \cite{chuamemristor71}. This 
characteristic was somehow lacking
in electrical circuit theory, since resistors, capacitors and inductors are defined
by voltage-current, charge-voltage and flux-current relations, respectively.
The design of such a {\em 
memory-resistor} or {\em memristor} at the nanometer scale
announced by an HP team in 2008 \cite{memristor2008}
has driven a lot of attention to these devices. 
Many applications are being developed
(cf.\  \cite{chen, itohchua09,  snider08, wang09, yang08} and references therein).
Memristive devices
pose challenging problems at the device and circuit modeling levels
\cite{benderli09,  pershin, memristors, wu}
but also from a dynamical perspective \cite{chuamemristor08, messias10,
muthus10, muthuskokate, membranch}.

The memristor reported in \cite{memristor2008} is governed by a relation
of the form
\begin{equation} \label{vi1}
v = M(q) i,
\end{equation}
being framed in the charge-control setting postulated by Chua in his
seminal paper \cite{chuamemristor71}; indeed, assuming $\phi$ to be $C^1$,
taking time derivatives in 
the equation $\varphi=\phi(q)$ and using the relations $\varphi'=v$, $q'=i$, we are led
to (\ref{vi1}) with $M(q)=\phi'(q)$.
Chua and Kang 
considered in \cite{chuakang} the more general characteristic
\begin{equation}
v = M(q, i) i. \label{vi2}
\end{equation}

Both (\ref{vi1}) and (\ref{vi2}), but also the maps governing other devices
such as capacitors or resistors, are
particular instances of the fully nonlinear relation
\begin{equation}
v = \eta(q, i) \label{vi3}
\end{equation}
to be discussed in this paper. 
The general characteristic (\ref{vi3}) will 
give a deeper insight into the mathematical properties that underly
several dynamical features of memristive circuits. In particular,
memristors increase the dynamical degree of freedom regardless
of their actual location in the circuit,
contrary to capacitors or inductors which in certain 
(so-called topologically degenerate)
configurations do {\em not}
increase the state dimension; additionally, memristors
have been observed to introduce 
null eigenvalues in the linearization of different 
circuits
\cite{messias10, membranch}. This kind of problems can be nicely
addressed in the framework here introduced.

This will be possible because of the fact that (\ref{vi3}), together with
its dual relation
\begin{equation}
i = \zeta(\varphi, v) \label{iv3}
\end{equation}
and their time-varying counterparts,
allow for 
a surprisingly simple, albeit general, 
model for the dynamics of
nonlinear circuits including 
resistors, 
memristors, 
capacitors,
inductors, and voltage and current sources. Indeed, 
such a model can be written as
\begin{subequations} \label{ceqs}
\begin{eqnarray}
q' & = & i_q \\
\varphi' & = & v_{\varphi} \\
v_q  & = & f(q, i_q, t) \label{ceqs3} \\ 
i_{\varphi} & = & g(\varphi, v_{\varphi}, t) \label{ceqs4} \\
0 & = & B_qv_q + B_\varphi v_{\varphi} \label{ceqs5} \\
0 & = & D_qi_q + D_\varphi i_{\varphi},\label{ceqs6} 
\end{eqnarray}
\end{subequations}  
where $B=(B_q \ B_\varphi)$ and $D=(D_q \ D_\varphi)$ stand for the 
reduced loop and cutset matrices of the
circuit (cf.\ subsection \ref{subsec-graphs}). 
The circuit elements (and hence $B$ and $D$) 
are divided into so-called $q$- and
$\varphi$-devices, as detailed in Section \ref{sec-mod}. The key aspect 
supporting the model (\ref{ceqs}) is that e.g.\ (\ref{ceqs4}) 
accommodates a great variety of devices, including 
flux-controlled memristors, linear
inductors, Josephson junctions, $pn$ and tunnel diodes, 
independent current sources, voltage-controlled current sources, etc.; 
similar remarks apply to (\ref{ceqs3}). 

The goal of our research is two-fold: first, from a 
modeling point of view, we want not only to accommodate in circuit theory
general devices governed by (\ref{vi3}) or (\ref{iv3}), but also 
emphasize the fact that these characteristics allow for a very general
description of the circuit dynamics, providing in turn a framework for
the analysis of different nonlinear phenomena; from an analytical standpoint, 
we wish to tackle several dynamical properties of circuits
with fully nonlinear memristors. The structure of the paper reflects the
two main aspects that drive this research: we address modeling issues in 
Section  \ref{sec-mod}, and dynamical 
properties in Section \ref{sec-ana}. Section \ref{sec-con} compiles some concluding
remarks.




\section{Device models}
\label{sec-mod}



\subsection{$q$-devices}

\begin{defin} A $q$-device is a circuit element which admits a $C^1$ description of the form
\begin{equation}
v = \eta(q, i, t). \label{qdev}
\end{equation}
\end{defin}

\vspace{1mm}

\noindent The terminology 
comes from the fact that $q$-devices may be
controlled by the charge $q$ (this will be the case for capacitors),
its time derivative $q'=i$ (for current-controlled resistors) 
or both (for $q$-memristors). 
Voltage sources will be included in this group for the
sake of completeness.

The relation (\ref{qdev}) is a time-varying generalization of (\ref{vi3}).
Depending on the actual form of the map $\eta$ in (\ref{qdev}), $q$-devices particularize to the ones
listed in Table \ref{table-q}. It is worth noting that $\eta$ need not be a scalar map; this accommodates coupling
effects within the different types of $q$-devices. 
Keep also in mind
that e.g. ${\partial \eta}/{\partial q} \not\equiv 0$
(resp.\ $\equiv 0$)
means that this derivative does not vanish (resp.\ does vanish) identically; in the first case,
it may find a zero at certain points, though; for instance, it vanishes at $i=0$
in Chua's memristor, for which $v=M(q)i$.

\vspace{5mm}
\begin{table}[htb]
\noindent \hspace{-2mm}
\begin{boxedminipage}{169mm}
\vspace{2mm}
\begin{enumerate}
\item 
A {\em $q$-memristor} 
is a $q$-device for which 
\begin{equation*} 
\frac{\partial \eta}{\partial q} \not\equiv 0, \ \ 
\frac{\partial \eta}{\partial i} \not\equiv 0.
\end{equation*}

\hspace{-7mm} The map governing 
the set of $q$-memristors will be denoted by
$v_m = \eta_1(q_m, i_m, t).$

\item A {\em capacitor} 
is a $q$-device for which 
\begin{equation*} 
\frac{\partial \eta}{\partial q} \not\equiv 0, \ \ 
\frac{\partial \eta}{\partial i} \equiv 0.
\end{equation*}
\hspace{-7mm} The relation governing the set of 
capacitors will be denoted by $v_c = \eta_2(q_c, t).$

\item A {\em current-controlled resistor} 
is a $q$-device for which 
\begin{equation*} 
\frac{\partial \eta}{\partial q} \equiv 0, \ \ 
\frac{\partial \eta}{\partial i} \not\equiv 0.
\end{equation*}
\hspace{-7mm} The map governing the set of 
current-controlled resistors will
be written as
$v_r = \eta_3(i_r, t).$

\item An independent {\em voltage source} 
is a $q$-device for which 
\begin{equation*} 
\frac{\partial \eta}{\partial q} \equiv 0, \ \ 
\frac{\partial \eta}{\partial i} \equiv 0.
\end{equation*}
\hspace{-7mm} The relation governing independent voltage sources will be denoted by
$v_u=\eta_4(t)$. \vspace{3mm}
\end{enumerate}
\end{boxedminipage}
\caption{$q$-devices.}\label{table-q}
\end{table}

In a time-invariant setting,
if $\eta$ is linear in $i$
then we get Chua's memristor, for which $v=M(q)i$. 
For this reason
we use the expression ``fully nonlinear'' to label memristors with the general
form $v=\eta(q, i)$. 
This should cause no misunderstanding with the nonlinear nature of the
original flux-charge relation $\varphi=\phi(q)$ in Chua's setting.
Note also that these memristors are referred to in the
literature
either as {\em charge-controlled memristors} or as 
{\em current-controlled memristors}; with the expression ``$q$-memristor''
we want to emphasize that the charge $q$ is the dynamic variable actually
involved in the description of the
device, as it will happen with the flux $\varphi$
in $\varphi$-memristors.

\subsection{$\varphi$-devices}

\begin{defin} A $\varphi$-device is a circuit element which admits a $C^1$ description of the form
\begin{equation}
i = \zeta(\varphi, v, t). \label{phidev}
\end{equation}
\end{defin}

\vspace{1mm}

\noindent 
The different types of $\varphi$-devices are enumerated 
in Table \ref{table-varphi}.

\begin{table}[htb]
\noindent \hspace{-2mm}
\begin{boxedminipage}{169mm}
\vspace{2mm}
\begin{enumerate}
\item A 
{\em $\varphi$-memristor}
is a $\varphi$-device for which 
\begin{equation*} 
\frac{\partial \zeta}{\partial \varphi} \not\equiv 0, \ \ 
\frac{\partial \zeta}{\partial v} \not\equiv 0.
\end{equation*}
\hspace{-7mm} The map governing the set of $\varphi$-memristors will be denoted by
$i_w = \zeta_1(\varphi_w, v_w, t).$

\item An {\em inductor} 
is a $\varphi$-device for which 
\begin{equation*} 
\frac{\partial \zeta}{\partial \varphi} \not\equiv 0, \ \ 
\frac{\partial \zeta}{\partial v} \equiv 0.
\end{equation*}
\hspace{-7mm} The relation governing the set of 
inductors will be denoted by
$i_l = \zeta_2(\varphi_l, t).$

\item A {\em voltage-controlled resistor} 
is a $\varphi$-device for which 
\begin{equation*} 
\frac{\partial \zeta}{\partial \varphi} \equiv 0, \ \ 
\frac{\partial \zeta}{\partial v} \not\equiv 0.
\end{equation*}
\hspace{-7mm} The map governing the set of 
voltage-controlled resistors will
be written as
$i_g = \zeta_3(v_g, t).$

\item An independent {\em current source} 
is a $\varphi$-device for which 
\begin{equation*} 
\frac{\partial \zeta}{\partial \varphi} \equiv 0, \ \ 
\frac{\partial \zeta}{\partial v} \equiv 0.
\end{equation*}
\hspace{-7mm} The relation governing independent current sources will be denoted by
$i_j=\zeta_4(t)$. \vspace{3mm}
\end{enumerate}
\end{boxedminipage}
\caption{$\varphi$-devices.}\label{table-varphi}
\end{table}

Again, Chua's flux-controlled memristor \cite{chuamemristor71}, for which $i=W(\varphi)v$,
is obtained in particular when $\zeta$ is time-independent and linear in $v$.

\



We do not impose any condition on the time derivatives, either for $q$-
or $\varphi$-devices. 
If they vanish in memristors, resistors, capacitors 
and inductors 
we would be led to a 
time-invariant setting, with $\eta_k$, $\zeta_k$ ($k=1, 2, 3$) being independent of
$t$.
For sources this assumption would model DC ones.
It is also worth remarking that
current-controlled voltage sources (CCVS's) could be easily included in the
above taxonomy as
$q$-devices and, similarly, voltage-controlled current sources (VCCS's) can be
modeled as $\varphi$-devices.

\subsection{Characteristic matrices and passivity}

The matrices of partial derivatives of the maps introduced above 
define the so-called  {\em characteristic matrices} of the different devices.

\begin{defin}
The incremental {\em memristance}, {\em elastance} and {\em
  resistance} matrices of 
$q$-memristors, capacitors and
current-controlled resistors, respectively, are defined as 
\begin{equation*}
M(q_m, i_m, t) =  \frac{\partial \eta_1(q_m, i_m, t)}{\partial i_m}, \ \
E_c(q_c, t) =  \frac{\partial \eta_2(q_c, t)}{\partial q_c}, \ \
R(i_r, t) =  \frac{\partial \eta_3(i_r, t)}{\partial i_r}. \ 
\end{equation*}
Additionally, the incremental {\em elastance} of $q$-memristors is
\[ 
E_m(q_m, i_m, t) =  \frac{\partial \eta_1(q_m, i_m, t)}{\partial q_m}.
\]\end{defin}

In Chua's memristor, for which $v=M(q)i$, the incremental
memristance depends only on the charge $q$. 
The fact that $q(t) = \int_{-\infty}^t i_m(\tau) d\tau$
makes this element behave as a resistor in which the resistance
depends on the device history; the ``memory-resistor'' name stems
from this.

If the elastance of capacitors is non-singular at a given operating point,
then because of the implicit function theorem these devices admit (at least
locally) a voltage-controlled description \[q_c=\gamma (v_c, t).\] The 
incremental {\em capacitance} matrix is then defined as
\begin{equation*}
C(v_c, t) =  \frac{\partial \gamma(v_c, t)}{\partial v_c}
=(E_c(\gamma (v_c, t), t))^{-1}.
\end{equation*}

Note that when a given set of devices (memristors, capacitors or resistors) does not exhibit
coupling effects, the corresponding characteristic matrix is a diagonal
one, with diagonal entries defined by the memristances, elastances or
resistances of the individual devices.

\

The characteristic matrices of $\varphi$-devices are defined analogously, as
detailed below.
\begin{defin}
The incremental {\em memductance}, {\em reluctance} and {\em
  conductance} matrices of 
$\varphi$-memristors, inductors and
voltage-controlled resistors, respectively, are defined as 
\begin{equation*}
W(\varphi_w, v_w, t) =  \frac{\partial \zeta_1(\varphi_w, v_w, t)}{\partial v_w}, \ \
{\mathcal R}_l(\varphi_l, t) =  \frac{\partial \zeta_2(\varphi_l, t)}{\partial \varphi_l}, \ \
G(v_g, t) =  \frac{\partial \zeta_3(v_g, t)}{\partial v_g}. \ 
\end{equation*}
The incremental {\em reluctance} of $\varphi$-memristors is
\[ 
{\mathcal R}_w(\varphi_w, v_w, t) =  \frac{\partial \zeta_1(\varphi_w, v_w, t)}{\partial \varphi_w}.
\]\end{defin}

As for capacitors, if the reluctance of inductors is non-singular,
then via the implicit function theorem these devices can be shown to
admit a local current-controlled description \[\varphi_l=\xi (i_l, t),\]
and the incremental {\em inductance} matrix is defined as
\begin{equation*}
L(i_l, t) =  \frac{\partial \xi(i_l, t)}{\partial i_l}
=({\mathcal R}_l(\xi (i_l, t), t))^{-1}.
\end{equation*}

In order to define the local notions of passivity and strict passivity, we make use
of the concept of a positive (semi)definite matrix; a square matrix $P$ is positive
semidefinite (resp.\ definite) if $u^T P u \geq 0$ (resp.\ $> 0$) for every
non-vanishing vector $u$. We do not assume $P$ to be symmetric.

\begin{defin}
We call the $q$-memristors, capacitors, current-controlled resistors,
$\varphi$-memris\-tors, inductors or voltage-controlled resistors {\em locally passive}
(resp.\ {\em strictly locally passive}) at a given point if the incremental 
memristance, elastance, resistance,
memductance, inductance or conductance matrix is positive semidefinite 
(resp.\ definite) at that point.
\end{defin}

It is easy to check that a positive definite matrix
is non-singular, and that its inverse is itself positive definite; this
means that the capacitance and inductance matrices are positive definite if the elastance
and reluctance matrices of capacitors and inductors are so.
No condition is imposed on the elastance and reluctance
matrices of
$q$- and $\varphi$-memristors, except for the fact that they cannot vanish
identically since otherwise these devices would amount to (current- or
voltage-controlled)
resistors. Note finally that, with terminological abuse, the adverb ``locally'' is often omitted
when a given set of devices satisfies the (semi)definiteness requirement
everywhere. Similarly, the term ``incremental'' is very often omitted in the 
characteristic matrices.

\section{Analytical results: nondegeneracy, null eigenvalues}
\label{sec-ana}

The general framework introduced in Section \ref{sec-mod} makes
it possible to address different analytical properties of
 memristive circuits in broad generality. After introducing some 
background material in subsection \ref{subsec-graphs},
we tackle in subsections \ref{subsec-complex},
\ref{subsec-null} and \ref{subsec-nullChua}
two particular problems which rely on linearization,
emphasizing the role of the form of the different 
characteristic matrices discussed above; specifically, 
the absence of memristance (resp.\ memductance)
in capacitors (resp.\ inductors), namely,
the conditions $\partial \eta / \partial i \equiv 0$,
$\partial \zeta / \partial v \equiv 0$ holding for them,
together with the form of the elastance and
reluctance matrices for different memristors, will explain
certain dynamical features of memristive circuits.
Our approach should also be useful in future analyses of 
many other aspects of memristive circuit dynamics, including 
e.g.\ stability properties, oscillatory phenomena or bifurcations.

\subsection{Some auxiliary results from digraph theory}
\label{subsec-graphs}

Many aspects of our analysis will rely on the properties of 
the {\em directed graph} or {\em digraph} underlying a given electrical
circuit. We compile below some background material on digraph theory.
The reader is referred to \cite{andras1, andras2, bollobas, foulds, wsbook}
for additional details.

We will work with a directed 
graph having $n$ nodes, $m$ branches and $k$ connected
components. A subset $K$ of the set of branches of a 
digraph is a {\em cutset} if the 
removal of $K$ increases the number of connected components of the digraph,
and it is minimal with respect to this property, that is,
the removal of any proper subset of $K$ does not increase the number of
components. In a connected digraph, a cutset is just 
a minimal disconnecting set of branches.
The removal of the branches of a cutset increases the number of connected
components by exactly one. 
Furthermore, all the branches of a cutset may be shown to connect the same
pair of connected components of the digraph which results from the deletion of
the cutset. 
This makes it possible to define the orientation of a cutset, say
from one of these components towards the other.
The cutset matrix $\tilde{D} = (d_{ij})$ is then defined by
\begin{eqnarray*}
d_{ij} = \left\{
\begin{array}{rl}
1 & \text{ if branch } j \text{ is in cutset } i \text{ with the same
  orientation }\\
-1 & \text{ if branch } j \text{ is in cutset } i \text{ with the opposite
  orientation } \\
0 & \text{ if branch } j \text{ is not in cutset } i.
\end{array}
\right.
\end{eqnarray*}
The rank of $\tilde{D}$ can be proved to be $n-k$; any set of 
$n-k$ linearly independent rows of $\tilde{D}$ defines
a {\em reduced cutset 
matrix} $D \in \R^{(n-k) \times m}$.
In a connected digraph, any reduced cutset matrix has order $(n-1) \times m$.

The space spanned by the rows of $D$ equals the one spanned by the rows of the so-called {\em incidence matrix}, defined as 
$A=(a_{ij})$ with
\begin{eqnarray*}
a_{ij} = \left\{
\begin{array}{rl}
1 & \text{ if branch } j \text{ leaves  node } i \\
-1 & \text{ if branch } j \text{ enters node } i \\
0 & \text{ if branch } j \text{ is not incident with node } i.
\end{array}
\right.
\end{eqnarray*}

Similarly, chosen an orientation in every loop, the 
{\em loop matrix} $\tilde{B}$ is defined as $(b_{ij})$, where
\begin{eqnarray*}
b_{ij} = \left\{
\begin{array}{rl}
1 & \text{ if branch } j \text{ is in loop } i \text{ with the same
  orientation }\\
-1 & \text{ if branch } j \text{ is in loop } i \text{ with the opposite
  orientation } \\
0 & \text{ if branch } j \text{ is not in loop } i.
\end{array}
\right.
\end{eqnarray*}
The rank of this matrix equals
$m-n+k$. A {\em reduced loop matrix} $B$
is any $((m-n+k)\times m)$-submatrix of $\tilde{B}$ with full row rank.

We denote by $B_K$ (resp.\ $B_{\mathcal{G}-K}$) the submatrix of $B$ 
defined by the columns which correspond to branches belonging 
(resp.\ not belonging) to a given set of branches
$K$;
the same applies to the
cutset matrix $D$. Certain submatrices of $B$ and $D$ 
characterize the existence of so-called $K$-cutsets
and $K$-loops (that is, cutsets or loops just defined by branches belonging
to  $K$), as stated below (cf.\ \cite[Sect.\ 5.1]{wsbook}).

\begin{lema} \label{lema-noloopscutsets}
Let $K$ be a subset of branches of a 
given digraph $\mathcal{G}$. The set $K$ does not contain cutsets
if and only if $B_{K}$ has full column rank or, equivalently,
iff $D_{\mathcal{G}-K}$ has full row rank.

Analogously, $K$ does not contain loops if and only if 
$D_K$ has full column rank or, equivalently,
iff $B_{\mathcal{G}-K}$ has full row rank.
\end{lema}

Actually, the dimensions of the spaces $\ke B_K$ and $\ke D_K$
are defined by the number of linearly independent $K$-cutsets and $K$-loops,
respectively.

The proof of the following result can be found e.g. in \cite[Sect.\ 7.4]{foulds}.
\begin{lema} \label{lema-orthogonalmatrices}
 If the columns of the reduced loop and
cutset matrices $B,$ $D$  of a digraph are 
arranged according to the same order of
branches, then $BD^T=0$, $DB^T=0$.
\end{lema}

The relations $\im D^T=\ke B$ and
$\im B^T=\ke D$ hold true. This expresses
that the {\em cut space}
$\im D^T$ spanned by the rows of $D$
can be described as $\ke B$ and, analogously,
the {\em cycle space} $\im B^T$ spanned by the rows of $B$ equals
$\ke D$ \cite{bollobas}. These spaces are orthogonal to each other
since $(\im D^T)^{\perp}= (\ke B)^{\perp}=\im B^T$. 

Finally, when applying these results to circuit analysis, we will
split the columns of $B$ and $D$ according to the nature of the devices 
accommodated in the corresponding branches; 
$B_q$, $B_{\varphi}$, $D_q$, $D_{\varphi}$ describe the submatrices of $B$
and $D$ defined by $q$- and
$\varphi$-devices; these matrices will be further split into
$B_m$, $B_c$, $B_r$, $B_u$, $B_w$, $B_l$, $B_g$, $B_j$ (resp.\
$D_m$, $D_c$ etc.), according to the division of 
 $q$- and
$\varphi$-devices into $q$-memristors, capacitors, current-controlled
resistors and voltage sources, and
$\varphi$-memristors, inductors, voltage-controlled
resistors and current sources, respectively. With this notation, the identity
$BD^T=0$ in Lemma \ref{lema-orthogonalmatrices} can be written as
\begin{eqnarray} \label{BD}
B_m D_m^T + B_c D_c^T + B_r D_r^T  + B_u D_u^T + B_w D_w^T + B_l D_l^T + B_g D_g^T + B_j D_j^T = 0.
\end{eqnarray}

\subsection{Nondegenerate problems and the order of complexity}
\label{subsec-complex}

The {\em order of complexity} or {\em state dimension}
of an electrical circuit is the number of
variables that can be freely assigned an initial value. 
Stemming from the work of
Bashkow and Bryant \cite{bashkow57, bryant59, bryant62}, 
a circuit composed of capacitors, inductors, resistors and
voltage and current sources is said to be {\em nondegenerate} if its 
order of complexity equals the number of reactive elements
(capacitors and inductors). A necessary condition for a circuit to be
nondegenerate
is its {\em topological} nondegeneracy, that is, the absence of
VC-loops (loops formed only by voltage sources and/or capacitors)
and IL-cutsets (cutsets defined only by current sources and/or inductors).
Different sufficient conditions can be given, involving e.g.\
the strict passivity of the circuit matrices or the structure of 
the circuit spanning trees \cite{lowindex, wsbook, normal, sommult, sommtree}.

In Theorem \ref{th-i1} 
we address the characterization of the order 
of complexity of circuits with fully nonlinear memristors, 
under strict passivity assumptions
and restricting the discussion to cases
without VC-loops and IL-cutsets. We show that memristors,
even under a fully nonlinear assumption,
do  {\em not} introduce topological degeneracies; this 
means that every memristor
increases by one the order of complexity, regardless of its 
location in the circuit.
The key difference with capacitors and inductors, which actually
introduce degeneracies when entering VC-loops or IL-cutsets,
is made by the fact that
the characteristics of these, namely, $v_c=\eta_2(q_c, t)$ and $i_l = \zeta_2(\varphi_l, t)$,
do not involve the current $i_c$ or the voltage $v_l$, respectively. 

Note, in this regard, that the equations $q_r'=i_r$, $q_u'=i_u$, $\varphi_g'=v_g$,
$\varphi_j'=v_j$ will be excluded from the dynamical description of
the circuit, since
the variables $q_r$, $q_u$, $\varphi_g$, $\varphi_j$ are decoupled from the rest
of the system and are hence irrelevant from the dynamical point of view. 
This means that the circuit dynamics will be defined by 
the model
\begin{subequations} \label{ceqsmod}
\begin{eqnarray}
q_{mc}' & = & i_{mc} \\
\varphi_{wl}' & = & v_{wl} \\
0 & = & v_q - f(q_{mc}, i_q, t) \label{ceqsmod3} \\ 
0 & = & i_{\varphi} - g(\varphi_{wl}, v_{\varphi}, t) \label{ceqsmod4} \\
0 & = & B_qv_q + B_\varphi v_{\varphi} \label{ceqsmod5} \\
0 & = & D_qi_q + D_\varphi i_{\varphi},\label{ceqsmod6} 
\end{eqnarray}
\end{subequations}  
in the
understanding that $q_{mc}$ and $\varphi_{wl}$ stand for ($q_m$, $q_c$) 
and ($\varphi_w$, $\varphi_l$), respectively. The same splitting applies
to $i_{mc}$ and $v_{wl}$. Notice that $f$ and $g$ group together
the maps $\eta_k$, $\zeta_k$ ($k=1, \ldots, 4$) arising in Tables
\ref{table-q} and \ref{table-varphi} for the sets of $q$-
and $\varphi$-devices.

\begin{teorema}\label{th-i1}
Consider a circuit without VC-loops and IL-cutsets in which
the memristance, resistance, memductance and conductance matrices $M$, 
$R$, $W$, $G$
are positive definite. Then its order of complexity is given by the total
number of $q$- and $\varphi$-memristors, capacitors, and inductors.
\end{teorema}

\noindent {\bf Proof.}
We will check that the relations (\ref{ceqsmod3})-(\ref{ceqsmod6})
make it possible to write explicitly all the branch voltages and currents in terms of
$q_m$, $q_c$, $\varphi_w$, $\varphi_l$, and that these relations impose no
constraint among these variables in the absence of VC-loops and
IL-cutsets. The way to do so is to check that the matrix of partial
derivatives of (\ref{ceqsmod3})-(\ref{ceqsmod6})
with respect to the branch voltages and currents $v_q$, $i_q$, $v_{\varphi}$, $i_{\varphi}$, 
is non-singular; a straightforward application of the implicit function
theorem yields the result.


This matrix of partial derivatives has the form
\begin{eqnarray} \label{i1matgen}
J = \left( \begin{array}{cccc}
I_q & -M_q & 0 & 0 \\
0 & 0 & -W_{\varphi}& I _{\varphi} \\
B_q & 0 & B_{\varphi} & 0 \\
0 & D_q & 0 & D_{\varphi} 
\end{array} \right)
\end{eqnarray}
with
\begin{eqnarray*} M_q =  \frac{\partial f}{\partial i_q}, \ \
W_{\varphi} = \frac{\partial g}{\partial v_{\varphi}}.
\end{eqnarray*}
Using a Schur reduction \cite{hor0, wsbook},
it is easy to see that (\ref{i1matgen}) is non-singular if and only if so it is
\begin{eqnarray} \label{i1matgenbis}
J_{\mathrm{red}}= \left( \begin{array}{cc}
B_q M_q & B_{\varphi} \\
D_q &  D_{\varphi}W_{\varphi} 
\end{array} \right).
\end{eqnarray}

According to the classification of $q$-devices into $q$-memristors, capacitors, resistors 
and voltage sources, the matrix $M_q$ has the block-diagonal structure block-diag$\{M, 0_c, R, 0_u\}$,
where $M$ and $R$ are the incremental memristance and resistance matrices. Analogously,
$W_{\varphi}$ reads as block-diag$\{W, 0_l, G, 0_j\}$. This confers the matrix
$J_{\mathrm{red}}$ in (\ref{i1matgenbis}) the form
\begin{eqnarray} \label{i1matbis}
\left( \begin{array}{cccccccc} 
B_m M  &   0 & B_rR  & 0 & B_w & B_l &  B_g & B_j \\
D_m  & D_c   & D_r  &  D_u & D_wW & 0 &  D_gG & 0 
\end{array}
\right).
\end{eqnarray}
Assume that
(\ref{i1matbis}) has  a non-trivial left-kernel; that is, suppose
that there exists a non-vanishing  $(x^T \ y^T)$ such that the following
relations hold:
\begin{subequations}
\begin{eqnarray}
x^T B_mM + y^T D_m & = & 0 \label{ker3} \\
y^T D_c & = & 0 \label{ker1} \\
x^T B_rR + y^T D_r & = & 0 \label{ker4} \\
y^T D_u & = & 0 \label{ker5}  \\
x^T B_w + y^T D_wW & = & 0 \label{kerbis} \\
x^T B_l & = & 0 \label{ker2} \\
x^T B_g + y^T D_gG & = & 0 \label{kertris} \\
x^T B_j & = & 0. \label{ker6} 
\end{eqnarray}
\end{subequations}
Multiply the identity (\ref{BD})
from the left by $x^T$ and from the right by $y$. 
Using (\ref{ker1}), (\ref{ker5}),
(\ref{ker2}) and (\ref{ker6}), 
we get
\begin{eqnarray*} 
 x^T(B_m D_m^T + B_r D_r^T + B_w D_w^T +  B_g D_g^T)y  = 0.
\end{eqnarray*} 
This equation can be recast,using (\ref{ker3}), (\ref{ker4}), (\ref{kerbis}) and
(\ref{kertris}),  as
\begin{eqnarray} \label{BDtris}
x^TB_m M^TB_m^Tx + x^TB_r R^TB_r^Tx  + y^T D_w W D_w^T y + y^T D_g G D_g^T y = 0.
\end{eqnarray} 
Be aware of 
the fact that $M^T$ and $R^T$ are positive definite since so they
are $M$ and $R$. Together with the positive definiteness of $W$ and $G$, 
this implies that  the relations 
\begin{subequations}
\begin{eqnarray}
x^T B_m & = & 0 \label{ker7} \\
x^T B_r & = & 0 \label{ker8} \\
y^T D_w & = & 0 \label{ker9} \\
y^T D_g & = & 0 \label{ker10} 
\end{eqnarray}
\end{subequations}
follow from (\ref{BDtris}). These relations yield, in turn,
\begin{subequations}
\begin{eqnarray}
y^T D_m & = & 0 \label{ker11} \\
y^T D_r & = & 0 \label{ker12} \\
x^T B_w & = & 0 \label{ker13} \\
x^T B_g & = & 0. \label{ker14} 
\end{eqnarray}
\end{subequations}
Now, the identities 
(\ref{ker2}),
(\ref{ker6}), (\ref{ker7}), (\ref{ker8}), (\ref{ker13}) and (\ref{ker14})
imply that $x=0$, because of  the absence of VC-loops and
Lemma \ref{lema-noloopscutsets}.
Similarly, 
(\ref{ker1}),
(\ref{ker5}), (\ref{ker9}), (\ref{ker10}), (\ref{ker11}) and (\ref{ker12}),
together with Lemma \ref{lema-noloopscutsets} and the absence of 
IL-cutsets, yield $y=0$. 
This means that (\ref{i1matbis}) and so $J$ in (\ref{i1matgen}) are
 indeed non-singular; therefore, the implicit function
theorem makes it possible to write $v_q$, $v_{\varphi}$, $i_q$, $i_{\varphi}$
in terms of $q_m$, $q_c$, $\varphi_w$ and $\varphi_l$, and this in turn yields
an explicit ODE formulated in terms of the charges of $q$-memristors and
capacitors
and the fluxes of $\varphi$-memristors and inductors. The state dimension of
this ODE is obviously defined by the number of memristors, capacitors and inductors and
the proof is complete.
\hfill $\Box$

\

Note that this result applies to fully nonlinear memristors,
regardless of the actual form of the maps $\eta_1$ and $\zeta_1$, as far as
the memristance 
$M$ and the memductance $W$ are positive definite. It
holds true, in particular, for Chua's memristors.

The different role played in this regard by $q$-memristors and capacitors,
on the one hand, and by $\varphi$-memristors and inductors, on the
other, becomes apparent in the light of the form of the matrix
(\ref{i1matbis}).
 Indeed, the
fact that the characteristic $v_c=\eta_2(q_c, t)$ of capacitors does not
involve the current $i_c$ 
is responsible for the vanishing block over $D_c$ in (\ref{i1matbis}).
Together with the null block over $D_u$, this explains the key role of
VC-loops in the state dimension problem. By contrast,
the presence of the $M$ matrix within (\ref{i1matbis})
makes the location of memristors irrelevant in this problem, showing that 
$q$-memristors cannot introduce topological degeneracies. In particular,
VCM-loops (with at least one $q$-memristor) do not reduce the state dimension
of the problem.
The same reasoning applies to $\varphi$-memristors, inductors,
IL-cutsets
and ILW-cutsets.

With more technical difficulties it is possible to show, assuming
the elastance of capacitors and the reluctance of inductors (or, equivalently,
the capacitance and the inductance) to be positive definite, that the order of 
complexity of topologically degenerate circuits with fully nonlinear
memristors is given by the total number
of memristors plus the number of capacitors in a normal tree 
(that is, a tree including all voltage sources, the maximum possible number of
capacitors, the minimum possible number of inductors and no current
source) and the number of inductors in a normal cotree. This provides
an  extension of Bryant's results \cite{bryant59, bryant62} to 
circuits with fully nonlinear memristors.

\subsection{Null eigenvalues; regular equilibrium points}
\label{subsec-null}

Many qualitative properties of dynamical systems at equilibria can be characterized in
terms of the linearized problem. In particular, zero eigenvalues of the
linearization may be responsible for stability changes and bifurcation
phenomena. In particular, within the context of memristive systems, some
particular circuits with Chua-type memristors have been shown in the literature to display null
eigenvalues \cite{messias10, membranch}. Below we tackle this problem for
general circuits with fully nonlinear memristors, addressing Chua's memristors as a
particular case 
with distinctive properties.

System (\ref{ceqsmod}) is a semiexplicit differential-algebraic
equation (DAE) \cite{bre1,kmbook,maeanew,rabrheintheo,wsbook}. Equilibrium
points are defined by the vanishing of its right-hand side, and the
linearization at a given equilibrium leads to the {\em matrix pencil}
\cite{gantma}
$\lambda H - K$, where $H=$ block-diag\{$I, 0$\} and $K$ is the 
Jacobian matrix of the right-hand side at equilibrium. The {\em spectrum} of the matrix pencil is the set of values of
$\lambda$ which make $\lambda H - K$ singular.
In particular, the pencil has null eigenvalues if and only if 
the Jacobian matrix $K$ is singular. 

An equilibrium point is said to be
{\em regular} if and only if $K$ is a non-singular matrix. Regular equilibria
are important regarding DC-solvability and Newton-based computations
\cite{dianafeldmann}. The following result
extends to systems with fully nonlinear memristors a property already known
for 
RLC circuits \cite{hagg86, hagg86b, matsuchua}. With the same terminological
convention, a VLW-loop is a loop defined by voltage sources, inductors and/or
$\varphi$-memristors, and an ICM-cutset is a cutset including only current
sources, capacitors and/or $q$-memristors.

\begin{teorema} \label{th-regular}
Consider a memristive circuit with positive definite incremental
resistance and conductance matrices $R$, $G$.
An equilibrium point of this circuit is regular
if and only if
\begin{itemize}
\item the elastances $E_m$, $E_c$ of $q$-memristors and capacitors, as well as the reluctances
  $\mathcal{R}_w$, 
$\mathcal{R}_l$ 
of $\varphi$-memristors and inductors, are non-singular at the equilibrium; and
\item the circuit does not have either VLW-loops or ICM-cutsets.
\end{itemize}
\end{teorema}

\noindent {\bf Proof.} The Jacobian matrix $K$ reads as
\begin{eqnarray} \label{Bmatrix}
K = \left( \begin{array}{cccccc}
0 & 0  & 0 & I_{mc} & I_{wl} & 0 \\
-E_{mc} & 0 & I_q & -M_q & 0 & 0 \\
0 & -\mathcal{R}_{wl} & 0 & 0 & -W_{\varphi}& I _{\varphi} \\
0 & 0 &B_q & 0 & B_{\varphi} & 0 \\
0 & 0 &0 & D_q & 0 & D_{\varphi} 
\end{array} \right)
\end{eqnarray}
with
\begin{eqnarray} \label{Emc} 
E_{mc} = \left( \begin{array}{cc}
E_m & 0 \\
0 & E_c \\
0 & 0 \\
0 & 0
\end{array} \right), \ 
\mathcal{R}_{wl} = \left( \begin{array}{cc}
\mathcal{R}_w & 0 \\
0 & \mathcal{R}_l \\
0 & 0 \\
0 & 0
\end{array} \right) 
\end{eqnarray}
and
\begin{eqnarray*} 
I_{mc} = \left( \begin{array}{cccc}
I_m & 0 & 0 & 0 \\
0 & I_c & 0 & 0 \\
0 & 0 & 0 & 0\\
0 & 0 & 0 & 0
\end{array} \right), \ 
I_{wl} = \left( \begin{array}{cccc}
0 & 0 & 0 & 0\\
0 & 0 & 0 & 0 \\
I_w & 0 & 0 & 0\\
0 & I_l & 0 & 0
\end{array} \right). 
\end{eqnarray*}

It is clear that 
the elastances $E_m$, $E_c$ and the reluctances ${\mathcal R}_w$, ${\mathcal
  R}_l$ must be non-singular for $K$
to be so; if these requirements are met, 
the non-singularity of $K$ relies on that of
\begin{eqnarray*} 
\left( \begin{array}{cccccccccc} 
0  &  0 & I_r &  -R  & 0 & 0 & 0 & 0 & 0 & 0  \\
0  & 0  & 0 &  0  & 0  &  -G  & 0 & 0 & 0 & I_g \\
B_m  &  B_c & B_r &  0   & 0 & B_g &  B_j & 0 & 0 & 0\\
0  &  0 &  0 &  D_r  &  D_u & 0 & 0 & D_w & D_l & D_g
\end{array}
\right)
\end{eqnarray*}
and, by means of a Schur reduction, on the non-singularity of
\begin{eqnarray} \label{Bredbis}
K_{\mathrm{red}}=\left( \begin{array}{cccccccc} 
B_m & B_c & B_rR  & 0 &  B_g & B_j & 0 & 0  \\
0 & 0   & D_r  & D_u  &  D_gG & 0 & D_w & D_l 
\end{array}
\right).
\end{eqnarray}
The blocks comprising $B_m$, $B_c$, $B_j$ show, in the light of Lemma 
\ref{lema-noloopscutsets},
that the absence of ICM-cutsets is necessary for $K_{\mathrm{red}}$ to be non-singular; analogously,
the blocks including $D_u$, $D_w$, $D_l$ rule out the
presence of VLW-loops. Conversely, in order to show that the absence of
these configurations guarantees $K_{\mathrm{red}}$ to be non-singular, it
suffices to check that a column-reordering gives this matrix the form of 
$J_{\mathrm{red}}$ in (\ref{i1matgenbis}); 
proceeding exactly as in the proof of Theorem \ref{th-i1}, and
using only the positive definiteness of $R$ and $G$, it is easy to check that
in the absence of VLW-loops and ICM-cutsets the matrix $K_{\mathrm{red}}$ (and
then $K$)
is actually non-singular, as we aimed to show.
\hfill $\Box$

\

More is true, actually; proceeding as in Theorem
\ref{th-zeroChua} 
below, one can show that, in problems in which 
the elastance of $q$-memristors and capacitors and the reluctance
of $\varphi$-memristors and inductors are non-singular, the
geometric multiplicity of the zero eigenvalue of $K$ equals the number
of independent VLW-loops and independent ICM-cutsets, where the notion of
independence relies on the cut- and cycle-spaces described in
subsection \ref{subsec-graphs}.

\subsection{Null eigenvalues in circuits with Chua-type memristors}
\label{subsec-nullChua}

A distinct feature of Chua's memristors is
that the elastance $E_m=\partial \eta_1 / \partial q_m$
and the reluctance $\mathcal{R}_w=\partial \zeta_1 / \partial \varphi_w$
do vanish at equilibrium points, because 
$\eta_1(q_m, i_m)=M(q_m)i_m$ and $\zeta_1(\varphi_w, v_w)=W(\varphi_w)v_w$ are
linear in $i_m$ and $v_w$, respectively, and the identities
$i_m=0$, $v_w=0$ hold at equilibria. Theorem \ref{th-regular} 
above then predicts the existence of 
null eigenvalues in Chua-type memristive circuits
owing to this zero-crossing property, as already observed in
specific examples \cite{messias10, membranch}. The number of memristive devices
actually characterizes the geometric multiplicity of the zero eigenvalue
in a broad class of circuits with Chua-type memristors,
as detailed below.

\begin{teorema} \label{th-zeroChua}
Consider a circuit in which the elastance $E_m$ of $q$-memristors and the
reluctance $\mathcal{R}_w$ of $\varphi$-memristors do vanish at equilibria.
If  the
resistance and conductance matrices $R$, $G$ are positive definite 
at a given equilibrium, and the elastance $E_c$ of capacitors and the 
reluctance $\mathcal{R}_l$ 
of inductors (equivalently, the capacitance and the inductance)
are non-singular, then the geometric multiplicity of 
the null eigenvalue equals the number of memristors,
independent VL-loops and independent IC-cutsets.
\end{teorema}

\noindent {\bf Proof.} Making $E_m=0$ and $\mathcal{R}_w=0$ in (\ref{Emc}),
the corank of the matrix $K$ can be easily seen to equal the number of
($q$- and $\varphi$-)
memristors plus the corank of 
\begin{eqnarray} \label{Bred3}
\left( \begin{array}{cccccccccccc} 
 0 & 0 & I_m  &  0   & 0 &  0  & 0 &  0 & 0 & 0 & 0 & 0  \\
 -E_c & 0 & 0  & I_c & 0 & 0 & 0 & 0 &   0 & 0 & 0 & 0  \\
 0 &  0 & 0  &  0 & I_r &  -R  & 0 & 0 & 0 & 0 & 0 & 0  \\
 0 &  0 & 0 & 0 & 0 & 0 &  0 & 0  & 0  & I_w & 0 & 0  \\
 0 &  -{\mathcal R}_l &0  & 0  & 0 &  0  & 0 &  0 & 0 & 0 & I_l & 0 \\
 0 &  0 &0  & 0  & 0 &  0  & 0  &  -G  & 0 & 0 & 0 & I_g \\
 0 &  0 &B_m  &  B_c & B_r &  0   & 0 & B_g &  B_j & 0 & 0 & 0\\
 0 &  0 &0  &  0 &  0 &  D_r  &  D_u & 0 & 0 & D_w & D_l & D_g
\end{array}
\right).
\end{eqnarray}
The corank of (\ref{Bred3}) in turn equals that of
\begin{eqnarray} \label{Bred4}
\left( \begin{array}{ccccccc} 
 B_c E_c &  0  & B_r R & 0 &  B_g & B_j \\
 0 &  D_l {\mathcal R}_l&  D_r & D_u  &  D_g G & 0 
\end{array}
\right).
\end{eqnarray}
Contrary to (\ref{Bredbis}), (\ref{Bred4}) is not a square matrix, and
its corank is defined by the dimension of its (right) kernel.
A vector $(u, v, w, x, y, z)$ of the kernel of this matrix verifies
\begin{subequations}\label{equs}
\begin{eqnarray}
B_c E_c u & = & 0 \\
B_r Rw & = & 0 \\
B_gy & = & 0 \\
 B_jz & = & 0   \\
D_l {\mathcal R}_l v & = & 0 \\
D_rw & = & 0 \\
 D_ux  & = & 0 \\  
D_g Gy & = & 0.
\end{eqnarray}
\end{subequations}
The relations (\ref{equs}) imply that
$(0, E_cu, Rw, 0, 0, 0, y, z) \in  \ke (B_m \ B_c \ B_r \ B_u \ B_w \ B_l
\ B_g \ B_j)= \ke B$, and
$(0, 0, w, x, 0, {\mathcal R}_l v, Gy, 0) \in \ke (D_m \ D_c \ D_r \ D_u \ D_w \ D_l
\ D_g \ D_j)=\ke D$.
Use then the identities $\ke B = \im D^T$ and 
$\ke D = \im B^T$ to derive the
existence of vectors $p$, $q$ satisfying
\begin{eqnarray*} 
0=D_m^T p, \ E_c u = D_c^T p, \ Rw = D_r^Tp, \ 0=D_u^T p, \ 0=D_w^T p, \ 0=D_l^T p, \ y = D_g^Tp, \ z=D_j^Tp,  \\
0=B_m^T q, \ 0 = B_c^Tq,  \  w= B_r^Tq, \ x=B_u^Tq, \ 0=B_w^Tq, \ {\mathcal R}_lv =
B_l^Tq, \ Gy = B_g^Tq, \ 0 = B_j^Tq. \
\end{eqnarray*}
Pre- and post-multiplying (\ref{BD}) by 
$q^T$ and $p$, respectively, we get,
in the light of these identities, 
\[w^T R w + y^T G^Ty =0,\] 
that is,
$w=y=0$ because of the positive definiteness of $R$ and $G$. Then
$(E_cu, z)$ belongs to 
$\ke (B_c \ B_j)$ and $(x, {\mathcal R}_lv)$ belongs to
$\ke (D_u \ D_l)$. Due to the non-singular nature of $E_c$ and ${\mathcal R}_l$,
this means that the corank of (\ref{Bred3}) equals the number
of independent IC-cutsets plus the number of independent VL-loops
(cf.\ subsection \ref{subsec-graphs}). 
Hence, the corank of the matrix
$K$, and therefore the geometric multiplicity of its null eigenvalue,
equals the number of ($q$- and $\varphi$-) memristors plus the 
number
of independent VL-loops and the number of independent IC-cutsets, as we
aimed to show.
\hfill $\Box$

\

In particular,
in circuits without VL-loops and IC-cutsets, the geometric
multiplicity of the null eigenvalue matches exactly the number of
memristors;
this is the case in the examples arising in \cite{messias10,membranch}.

\section{Concluding remarks}
\label{sec-con}

The fully nonlinear characteristics (\ref{vi3}) and (\ref{iv3}) introduced
in the present paper might accommodate future 
devices arising in nonlinear circuit theory and displaying
memristive effects, or provide more accurate models for already
existing devices. The general form of
these relations and their time-varying counterparts allows for a 
simple description of a very broad class of electrical and electronic
circuits; different dynamical features of these can be then addressed
in great generality. Several of these features can be understood
to follow from the special form of the characteristic matrices 
of the different devices arising in the analysis.
For instance, the zero-crossing property of Chua-type memristors 
makes them display a vanishing elastance and reluctance at equilibria;
from Theorems \ref{th-i1} and \ref{th-zeroChua},
it then follows that
in strictly passive circuits with Chua-type memristors,
not exhibiting VC-loops, IL-cutsets, VL-loops or
IC-cutsets, the order of complexity equals the number
of memristors plus reactive elements, every memristor introducing
a vanishing natural frequency in the linearized problem. Similarly, 
the absence of null eigenvalues is characterized, for general memristive
circuits and in terms of the circuit matrices and the
digraph topology, in Theorem \ref{th-regular}. The characterization 
of other analytical properties of memristive circuits within the framework
here introduced defines the scope of future research.

\end{document}